\documentclass[pre,aps,twocolumn,showpacs,floatfix]{revtex4-1}

\usepackage{amssymb}
\usepackage{graphicx}

\begin{document}

\title{Generalized chronotaxic systems: time-dependent oscillatory dynamics stable under continuous perturbation}

\author{Yevhen F. \surname{Suprunenko}}
\author{Aneta \surname{Stefanovska}}
\affiliation{Department of Physics, Lancaster University, Lancaster LA1 4YB, United Kingdom}
\date{\today}

\begin{abstract}
Chronotaxic systems represent deterministic nonautonomous oscillatory systems which are capable of resisting continuous external perturbations while having a complex time-dependent dynamics. Until their recent introduction in \emph{Phys. Rev. Lett.} \textbf{111}, 024101 (2013) chronotaxic systems had often been treated as stochastic, inappropriately, and the deterministic component had been ignored.
While the previous work addressed the case of the decoupled amplitude and phase dynamics, in this paper we develop a generalized theory of chronotaxic systems where such decoupling is not required.
The theory presented is based on the concept of a time-dependent point attractor or a driven steady state and on the contraction theory of dynamical systems.
This simplifies the analysis of chronotaxic systems and makes possible the identification of chronotaxic systems with time-varying parameters.
All types of chronotaxic dynamics are classified and their properties are discussed using the nonautonomous Poincar\'e oscillator as an example.
We demonstrate that these types differ in their transient dynamics towards a driven steady state and according to their response to external perturbations.
Various possible realizations of chronotaxic systems are discussed, including systems with temporal chronotaxicity and interacting chronotaxic systems.
\end{abstract}


\maketitle

\section{Introduction}

Complex dynamics, observed in real physical systems, is often modeled as stochastic or chaotic, or high-dimensional autonomous. Such a description is inappropriate when these systems are open (which is almost always the case), i.e. when they depend on time explicitly or when they are exposed to continuous perturbation originating from the external environment which cannot be considered as a part of a system.
In such cases a system should be described and modeled as nonautonomous. Nonautonomous dynamics is common in nature. The necessity to understand such a dynamics stimulated developments in the theory of nonautonomous dynamical systems \cite{Kloeden:11,Vishik:92} and the closely related theory of random dynamical systems \cite{Crauel:94,Romeiras:90}, as well as various inverse approach methods  \cite{Daubechies:11,Jamsek:10,Sheppard:11,Sheppard:12a,Stankovski:12b,Clemson:14a}. Thus, nonautonomous dynamics was studied in the life sciences \cite{Kloeden:13book}, in neuroscience \cite{Stam:05,Friston:12}, in cardiovascular \cite{Shiogai:10} and cardiorespiratory systems \cite{Iatsenko:13}, in cells \cite{Kurz:10}, in climate \cite{Bretherton:99,Chekroun:11,Ashwin:12}, and in solid state physics \cite{Konstantinov:11}.
Moreover, the description of a system as nonautonomous is the only description which is capable of explaining the \emph{stability of time-dependent dynamics}, -- such stability does not allow the time-variable dynamics to be changed easily by external \emph{continuous} perturbation. We will show that such stability allows the nonautonomous dynamics to look stochastic-like and vary complex, which may cause their misidentification as stochastic or chaotic. However, the underling deterministic dynamics of such systems is kept stable. Also, such properties allow the dynamics to be decomposable into deterministic underlying dynamics and dynamics due to external perturbations.

\begin{figure}[]
 \includegraphics[width=\columnwidth]{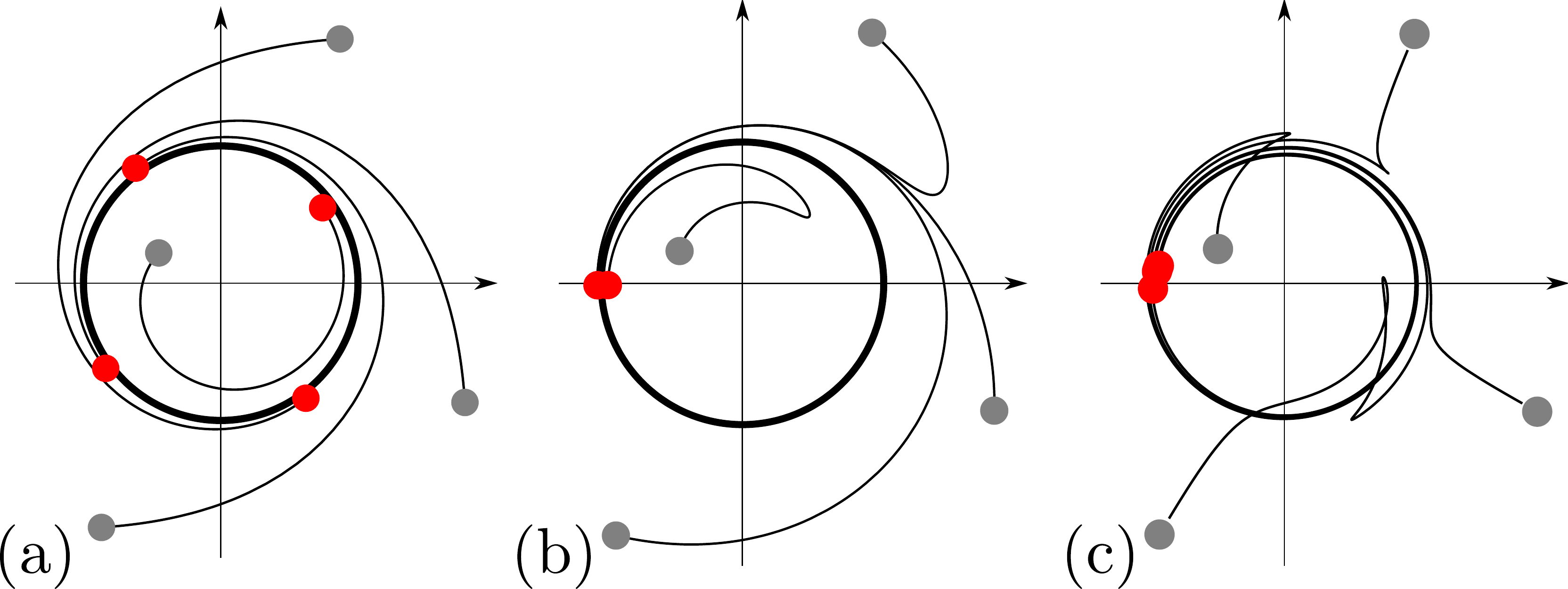}
\caption{(color online) (a) Conventional limit cycle (black circle). Oscillations with different initial conditions (gray points) eventually have the same amplitude but different phases. (b) Chronotaxic limit cycle (black circle). Initial conditions are completely forgotten. Amplitude and phase dynamics are decoupled. (c) Chronotaxic systems introduced in this paper have coupled phase and amplitude dynamics.
}
\label{fig: 0}
\end{figure}

Special attention has been drawn to time-dependent oscillatory dynamics in living systems \cite{Glass:01,Stefanovska:07} which can self-adjust and self-organise \cite{Friston:12}, and their time-dependent oscillatory dynamics can resist continuous external perturbations. Describing oscillations by their amplitude and phase, observations in living \cite{Shiogai:10,Kurz:10} and in solid state \cite{Konstantinov:11} systems show features of stability not only in their time-dependent amplitude dynamics, but also in their phase dynamics. Such stability in a phase dynamics cannot be explained by conventional autonomous models, such as limit-cycle models of self-sustained oscillations \cite{Andronov:66,Kuramoto:84,Pikovsky:01} -- a limit-cycle model provides stability only to the amplitude but not to the phase. In a limit cycle, the phase can easily be shifted by external perturbations, and such a shift does not decay. Such a typical dynamics is shown in Fig.~\ref{fig: 0}(a), where the amplitude and phase of the oscillations correspond to polar coordinates, and a limit cycle is a circle.

Assuming a decoupling of the time-dependent amplitude and phase dynamics, stability in the dynamics of each was described by the chronotaxic system model recently introduced in Refs.~\cite{Suprunenko:13,Suprunenko:14}. Chronotaxic systems represent a sub-class of nonautonomous oscillatory systems. Their main property is a stability in phase and amplitude, which leads to a decay of perturbations and, consequently, initial conditions become unimportant. Such a property emerges because of a time-dependent point attractor or a driven steady state in both the phase and amplitude dynamics. This results in a chronotaxic limit cycle, its typical dynamics is shown in Fig.~\ref{fig: 0}(b) -- the phase and amplitude independently approach their time-dependent point attractors.

In this paper we extend the theory of chronotaxic systems: we develop a mathematical definition of high-dimensional chronotaxic systems where the amplitude and phase dynamics do not need to be separated, in contrast to previous work \cite{Suprunenko:13,Suprunenko:14}. One example of such chronotaxic dynamics is shown in Fig.~\ref{fig: 0}(c), where radial and angular transient dynamics influence each other. The definition of chronotaxic systems presented in this paper is based on the concept of a time-dependent point attractor \cite{Kloeden:11} and on the contraction theory of dynamical systems \cite{Lohmiller:98,Pham:07,Sontag:10}. The use of a contraction theory, demonstrated on the example of a nonautonomous Poincar\'e oscillator, simplifies the analysis and identification of chronotaxic systems, and it makes possible the identification of chronotaxic systems with time-varying parameters possible. First, using constant parameters we build a diagram which denotes regions in a parameter space where the system is chronotaxic. Then, considering time-dependent parameters and analysing the dynamics of a contraction region, we find how parameters can depend on time to preserve chronotaxicity of the system.
A rich variety of types of generalized chronotaxic dynamics is found. This variety of types is due to the different combinations of (a) a point attractor with (b) a nonautonomous analogue of a limit cycle and (c) regions where trajectories diverge. Using examples we demonstrate that these types differ in their response to external perturbations. Different realizations of chronotaxic systems are discussed.

The theory of chronotaxic systems, presented in this paper and in Refs.~\cite{Suprunenko:13,Suprunenko:14}, together with corresponding inverse approach methods \cite{Clemson:14b}, make it possible to identify the underlying deterministic dynamics within the complex stochastic-like dynamics. It will be useful in various research fields, especially in living systems, where the identification of systems with complex stochastic-like dynamics as chronotaxic can help us in understanding their structure and function and their interactions with the external environment.

The paper is structured as follows: Sec.~II provides a generalized theory of chronotaxic systems, and Sec.~III provides an analysis of nonautonomous Poincar\'e oscillator as an example of a chronotaxic system with non-separable amplitude and phase dynamics.
In Sec.~IV a realization of chronotaxic systems is discussed. A summary of results is shown in Sec.~V.

\section{Generalized theory of chronotaxic systems}
A chronotaxic system \cite{Suprunenko:13,Suprunenko:14} is a nonautonomous oscillatory dynamical system $\mathbf{x}$ generated by an autonomous system of unidirectionally coupled equations
\begin{equation}
 \dot{\mathbf{p}}=\mathbf{f}(\mathbf{p});\qquad
 \dot{\mathbf{x}}=\mathbf{g}(\mathbf{x},\mathbf{p}).
 \label{1}
\end{equation}
Here $\mathbf{p}\in R^{n},$ $\mathbf{x}\in R^m,$ $\mathbf{f}:R^n\rightarrow R^n,$ $\mathbf{g}:R^m\rightarrow R^m.$ The system (\ref{1}) may also be called a master-slave configuration \cite{Haken:04}, or a drive and response system \cite{Kocarev:96}. Within the theory of nonautonomous dynamical systems \cite{Kloeden:01,Kloeden:11,Kloeden:12,Kloeden:12a} the system (\ref{1}) can be viewed as a skew-product flow or as a process. The solution $\mathbf{x}(t,t_0,\mathbf{x}_0)$ of Eqs.~(\ref{1}), which can be viewed as nonautonomous in the sense that $\dot{\mathbf{x}}=\mathbf{g}(\mathbf{x},\mathbf{p}(t)),$ depends on the actual time $t$ as well as on the initial conditions $(t_0,\mathbf{x}_0).$ For all $(t,t_0,\mathbf{x}_0)\in R\times R\times R^m$ the solutions satisfy the initial conditions, $\mathbf{x}(t_0,t_0,\mathbf{x}_0)=\mathbf{x}_0,$ and the co-cycle property, $\mathbf{x}(t_2,t_0,\mathbf{x}_0)=\mathbf{x}(t_2,t_1,\mathbf{x}(t_1,t_0,\mathbf{x}_0)).$

Following Ref.~\cite{Suprunenko:13}, the defining property of a chronotaxic system is the stability of its time-dependent dynamics in face of external perturbations. It is realized by a single time-dependent steady state or point attractor which performs an oscillatory motion. Following the theory of nonautonomous dynamical systems \cite{Kloeden:11}, a time-dependent steady state or point attractor is defined as a point $\mathbf{x}^{A}(t)$ in a state space, and it satisfies the following conditions of forward and pullback attractions,
\begin{eqnarray}
 \lim_{t\rightarrow+\infty}|\mathbf{x}(t,t_0,\mathbf{x}_0)-\mathbf{x}^{A}(t)|=0;\label{2}\\
 \lim_{t_0\rightarrow-\infty}|\mathbf{x}(t,t_0,\mathbf{x}_0)-\mathbf{x}^{A}(t)|=0,\label{3}
\end{eqnarray}
and a condition of invariance
\begin{equation}
 \mathbf{x}(t,t_0,\mathbf{x}^{A}(t_0))=\mathbf{x}^{A}(t).
 \label{condition invariance}
\end{equation}
Thus, the attraction of a system's state towards $\mathbf{x}^{A}(t)$ means that the initial conditions are forgotten. According to Eq.~(\ref{condition invariance}), the time-dependent steady state $\mathbf{x}^{A}(t)$ should be a solution of the system.

\begin{figure}[]
 \includegraphics[width=\columnwidth]{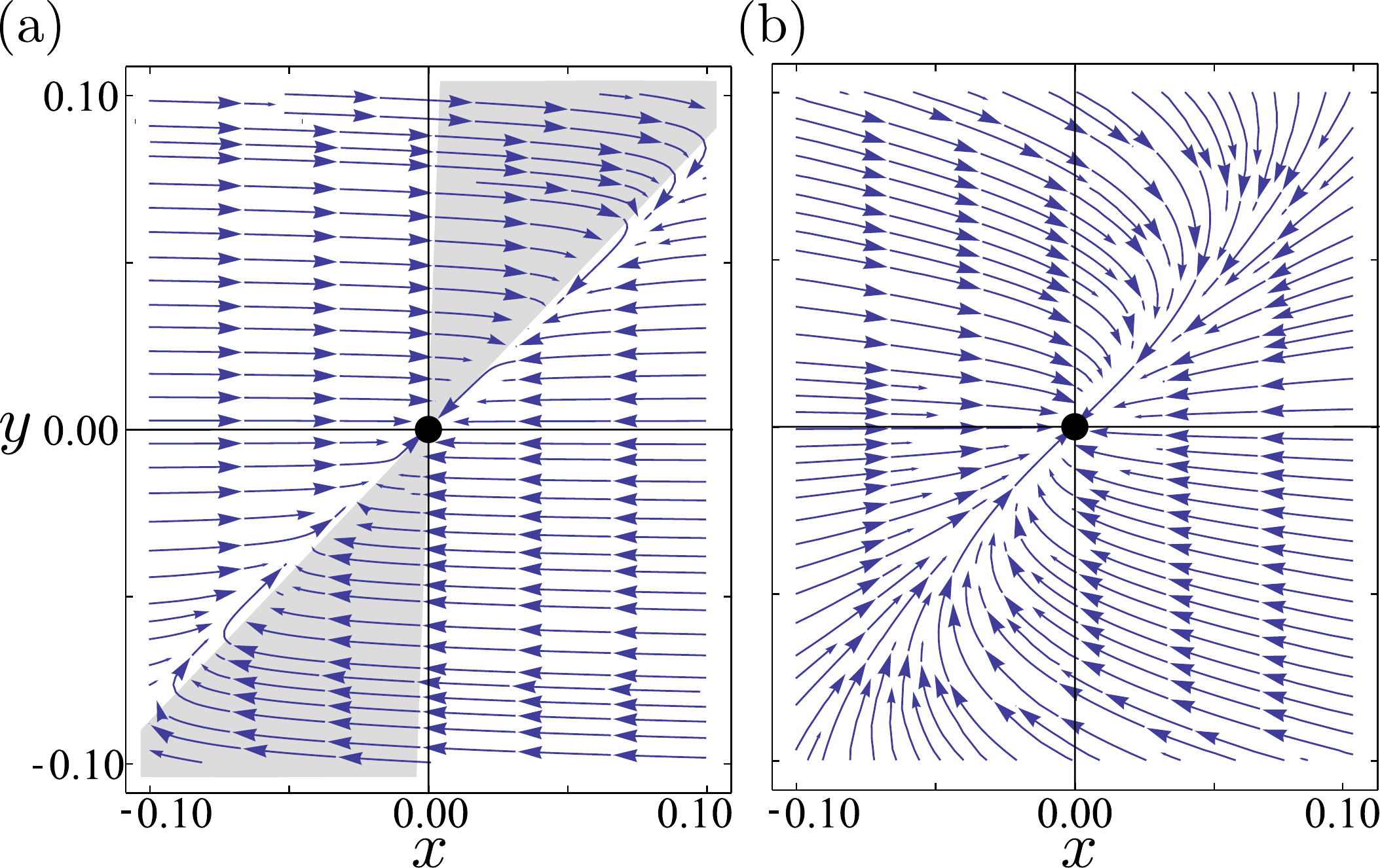}
\caption{Phase portraits of systems where a distance to a stationary fixed point (black dot) can: (a) temporarily increase (the area of the phase space where it occurs is gray); (b) only decrease.}
\label{fig: dev}
\end{figure}

It is important to stress that, in chronotaxic systems, a time-dependent steady state $\mathbf{x}^{A}(t)$ attracts other states of the system at all times. This means that in the unperturbed chronotaxic system the infinitesimal deviations from a steady state can only decrease in time. This condition was defined in systems with separable dynamics \cite{Suprunenko:13,Suprunenko:14} by linear stability analysis of a point attractor. However, this is not sufficient in the case of two and higher dimensional chronotaxic systems with coupled amplitude and phase. This can be demonstrated by the simple linear system $\dot{x}=-4 x+ 4.75 y;$ $\dot{y}=-0.2y$, whose phase portrait is shown in Fig.~\ref{fig: dev}(a). Its point attractor (fixed point $(x,y)=(0,0))$ is stable, as the eigenvalues are both negative, $-0.2$ and $-4$.
Nevertheless, when the system is in the gray area in Fig.~\ref{fig: dev}(a), its distance to the point attractor increases temporarily. Thus, strictly speaking, such a system does not resist perturbations, as some deviations due to perturbations can temporarily increase. This does not happen, however, in the example shown in Fig.~\ref{fig: dev}(b), which corresponds to $\dot{x}=-4 x+ 3.125 y;$ $\dot{y}=-1.5y$. The eigenvalues here are $-1.5$ and $-4.$ Thus, the negativity of characteristic exponents is not in itself enough to identify a system where any perturbation will continuously decay.

In order to take this into account, we use a contraction analysis instead of a linear stability analysis. Using a small virtual displacement $\delta \mathbf{x}^{A}(t)$ from the point attractor $\mathbf{x}^{A}(t)$, $\delta \mathbf{x}^{A}(t)=\mathbf{x}(t)-\mathbf{x}^{A}(t),$ one can write
\begin{equation}
 {d\over dt}\delta \mathbf{x}^A(t)={\partial \mathbf{g}(\mathbf{x}^A,t)\over \partial \mathbf{x}}\delta \mathbf{x}^A(t).
\end{equation}
The square of the distance $|\delta\mathbf{x}^A|^2=(\delta \mathbf{x}^A)^T\delta\mathbf{x}^A$ should decay in time, thus
\begin{equation}
 {d\over dt}|\delta\mathbf{x}^A|^2
 =
 2\delta \mathbf{x}^{A~T} J(\mathbf{x}^A,t)\delta\mathbf{x}
 <0.
\label{condition contraction}
\end{equation}
Here
\begin{equation}
J(\mathbf{x},t)={1\over 2}\left(\partial \mathbf{g}/ \partial \mathbf{x}+\partial \mathbf{g}/ \partial \mathbf{x}^{T}\right).
\label{Jacobian}
\end{equation}
 Denoting the maximum eigenvalue of the Jacobian $J$ by $\lambda_{max}(\mathbf{x}^A,t),$ the stability of $\mathbf{x}^{A}(t)$ at all times is constituted by the requirement that $\lambda_{max}$ is uniformly strictly negative, i.e.
\begin{equation}
 \exists \beta>0,~\forall t\geq 0,~
 \lambda_{max}(\mathbf{x}^{A}(t),t)\leq-\beta <0.
 \label{condition contraction2}
\end{equation}
This is the condition of contraction (similar to the requirement of having negative exponents in a linear stability analysis as in Ref.~\cite{Suprunenko:14}). As a result, any infinitesimal deviation from the point attractor will decay exponentially. The conditions (\ref{condition contraction2}) and (\ref{2})-(\ref{condition invariance}) determine chronotaxic systems.

The definition of chronotaxic systems is based on the explicit notion of a point attractor, and is difficult to use when data (numerically generated or experimental) are available. All of the conditions presented require prior knowledge of a point attractor. To avoid this inconvenience one needs to analyze stability in a stronger sense. Such an approach to stability is provided by contraction theory \cite{Lohmiller:98,Pham:07} and the theory of incremental stability \cite{Angeli:02,Ruffer:13}. These theories study dynamical systems where the distance between any two trajectories decays in time. Thus instead of a stable point it is important to introduce a region of contraction, which is defined as follows \cite{Lohmiller:98,Pham:07,Russo:09}. The region $C$ is called a contraction region for a system $\dot{\mathbf{x}}=\mathbf{g}(\mathbf{x},t)$ if
\begin{eqnarray}
 &&\forall \mathbf{x}\in C~~\exists\beta>0,\forall t~:\nonumber \\
 &&{1\over 2}\left({\partial \mathbf{g}(\mathbf{x},t)\over \partial \mathbf{x}}+
{\partial \mathbf{g}(\mathbf{x},t)\over \partial \mathbf{x}}^{T}\right)\leq-\beta I<0.
\label{contr region}
\end{eqnarray}
In contracting \cite{Lohmiller:98,Pham:07} and incrementally stable \cite{Angeli:02,Ruffer:13} systems, where a contraction region occupies the whole space, there is a unique time-dependent point attractor. This follows from conditions (\ref{2})-(\ref{condition invariance}) and (\ref{condition contraction2}).

\begin{figure}[]
 \includegraphics[width=0.9\columnwidth]{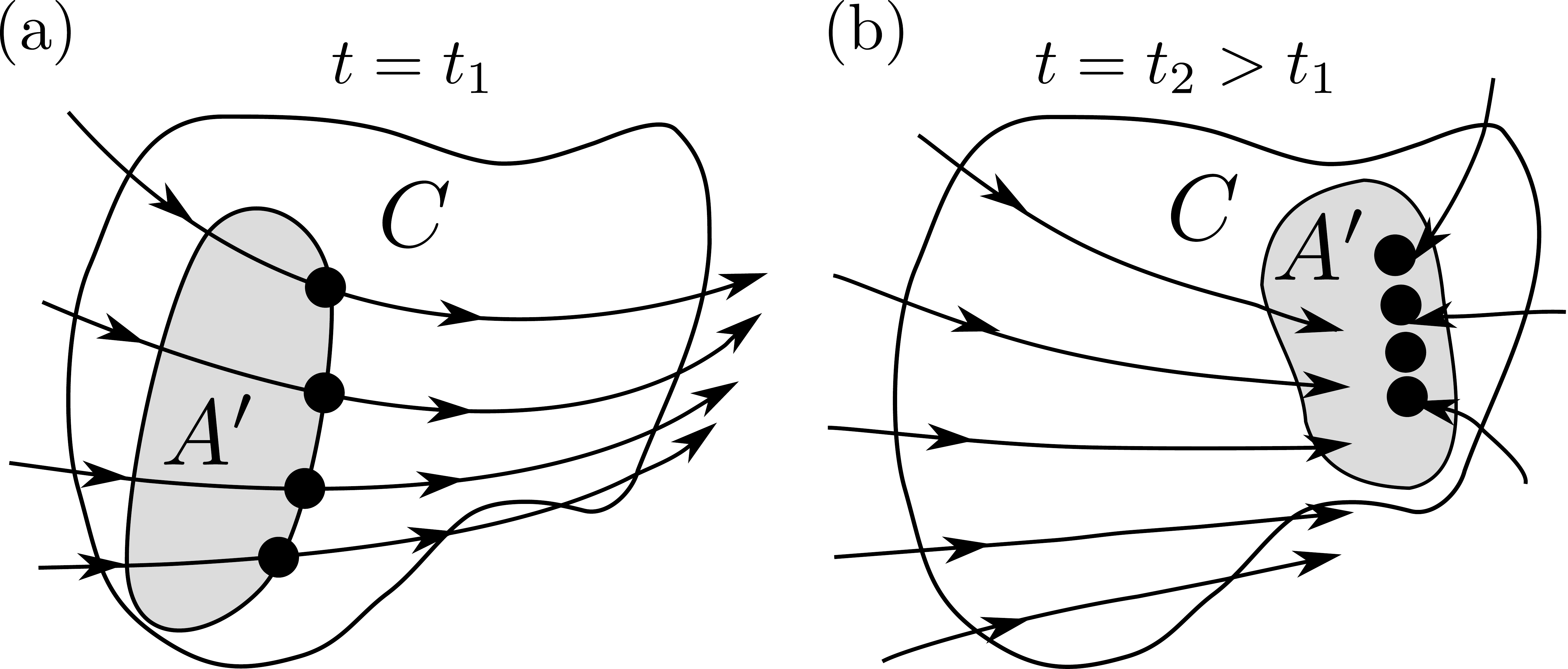}
\caption{Sketch of a contraction region $C$ in the phase space of a chronotaxic system at (a) time $t_1,$ (b) time $t_2>t_1.$ Arrows denote the velocity field. States of a system (white points) cannot leave an area $A^{\prime}$ which lies inside the contraction region $C.$ A time-dependent point attractor exists inside $A^{\prime}.$}
\label{fig: C}
\end{figure}

In chronotaxic systems the contraction region can be finite and can move in phase space. However, the fulfillment of the conditions (\ref{2})-(\ref{condition invariance}) and (\ref{condition contraction2}) imposes certain restrictions. Thus, inside the contraction regions there should be a finite area $A^{\prime}\subset C$ consisting of solutions of the system which never leave this finite area $A^{\prime},$ see Fig.~\ref{fig: C}. This requirement leads to the existence of a point attractor which satisfies the conditions (\ref{2})-(\ref{condition invariance}) inside the area $A^{\prime}\subset C.$  As it is a contraction regions, distances between trajectories inside this region can only decrease. Consequently, such a system will resist continuous external perturbation, and it will keep its time-dependent dynamics stable.

The definition of chronotaxic systems can therefore be reformulated without prior knowledge of a point attractor. \emph{A chronotaxic system is an oscillatory nonautonomous dynamical system which has a contraction region $C$ determined by Eq.~(\ref{contr region}) which contains a finite non-zero area $A^{\prime}\subset C$ such that states of a system inside $A^{\prime}$ do not leave $A^{\prime}.$} As discussed below, this definition is useful for the identification of chronotaxic systems from observations or from corresponding dynamical equations. It will be demonstrated and discussed in the next section.

Considering the relationship to known classes of dynamical systems, one may note that chronotaxic systems represent a subclass of asymptotically stable \cite{He:92,Kocarev:96} and converging \cite{Pavlov:04,Ruffer:13} systems. From the definition of chronotaxic systems it follows that chronotaxic systems are asymptotically stable (i.e. $\forall \mathbf{x}_{01,} \mathbf{x}_{02}:~ \lim_{t\rightarrow\infty}||\mathbf{x}(t,t_0,\mathbf{x}_{01})-\mathbf{x}(t,t_0,\mathbf{x}_{02})||=0$).
However, the requirement to have a point attractor at all times (before the limit $t\rightarrow\infty$ is achieved) is not necessarily fulfilled in asymptotically stable systems.
Additional conditions should also be satisfied, which is the existence of a pullback attractor (\ref{3}), its invariance (\ref{condition invariance}), and contraction (\ref{condition contraction2}). Thus, general asymptotically stable systems can temporarily become divergent, i.e. the contraction region and point attractor can disappear. Therefore, chronotaxic systems are a subclass of asymptotically stable systems.
Converging systems are defined as asymptotically stable systems with bound solutions \cite{Pavlov:04,Ruffer:13}. Therefore chronotaxic systems with bound solutions are also a subclass of converging systems.
Contracting \cite{Lohmiller:98,Pham:07} and incrementally stable \cite{Angeli:02,Ruffer:13} systems are most similar to chronotaxic systems as they satisfy contraction condition in a stronger sense discussed above, i.e. the distance between any two states of the system can only decrease in time. However, in contrast to contracting or incrementally stable systems, chronotaxic systems can have non-contracting regions in phase space, and these regions can move. Therefore, the distance between some points can increase temporarily. As a result, chronotaxic systems can generate a complex dynamics when perturbed.

\section{Nonautonomous Poincar\'e oscillator as a chronotaxic system}
Properties and different types of chronotaxic dynamics can be demonstrated using a simple but general example of an oscillatory self-sustained system. Consider a Poincar\'e oscillator that is unidirectionally coupled with a coupling strength $\varepsilon_{_{A}}(t)$ to a moving point with Cartesian coordinates $x_p=r_p\cos\alpha_p(t)$ and $y_p=r_p\sin\alpha_p(t)$:
\begin{eqnarray}
  \dot{x}&=&\varepsilon_{_{\Gamma}}\left(r_{p}-\sqrt{x^2+y^2}\right)x-\omega_0 y
  -\varepsilon_{_{A}}(t)\left(x-r_{p}\cos\alpha_p(t)\right),
\nonumber\\
  \dot{y}&=&\varepsilon_{_{\Gamma}}\left(r_p-\sqrt{x^2+y^2}\right)y+\omega_0 x
  -\varepsilon_{_{A}}(t)\left(y-r_{p}\sin\alpha_p(t)\right).
\nonumber\\
 \label{general example}
\end{eqnarray}
Here $\omega_0$ is a natural frequency, and $\varepsilon_{_{\Gamma}}$ is a parameter. When dealing with such a nonautonomous oscillatory system, it is important to find when the system is chronotaxic, i.e. to find the corresponding restrictions on the time-dependent coupling strength $\varepsilon_{_{A}}(t)$ and the frequency $\omega_p(t),$
\begin{equation}
 \omega_p(t)\equiv {d\alpha_{p}(t)\over dt}.
\end{equation}
The analytic derivation of the chronotaxicity conditions for such systems is difficult and in most cases is impossible. Here we show how the definition of chronotaxic systems presented in this paper allows us to find a region in parameter space where the system (\ref{general example}) is chronotaxic. Moreover, it becomes possible to find restrictions on the time-dependent functions $\varepsilon_{_{A}}(t)$ and $\omega_p(t)$ such that the system is chronotaxic. For this, first the application of contraction theory to a simple case without a coupling $\varepsilon_{_{A}}(t)$ is demonstrated.

\subsection{The case without coupling}
\begin{figure}[]
 \includegraphics[width=\columnwidth]{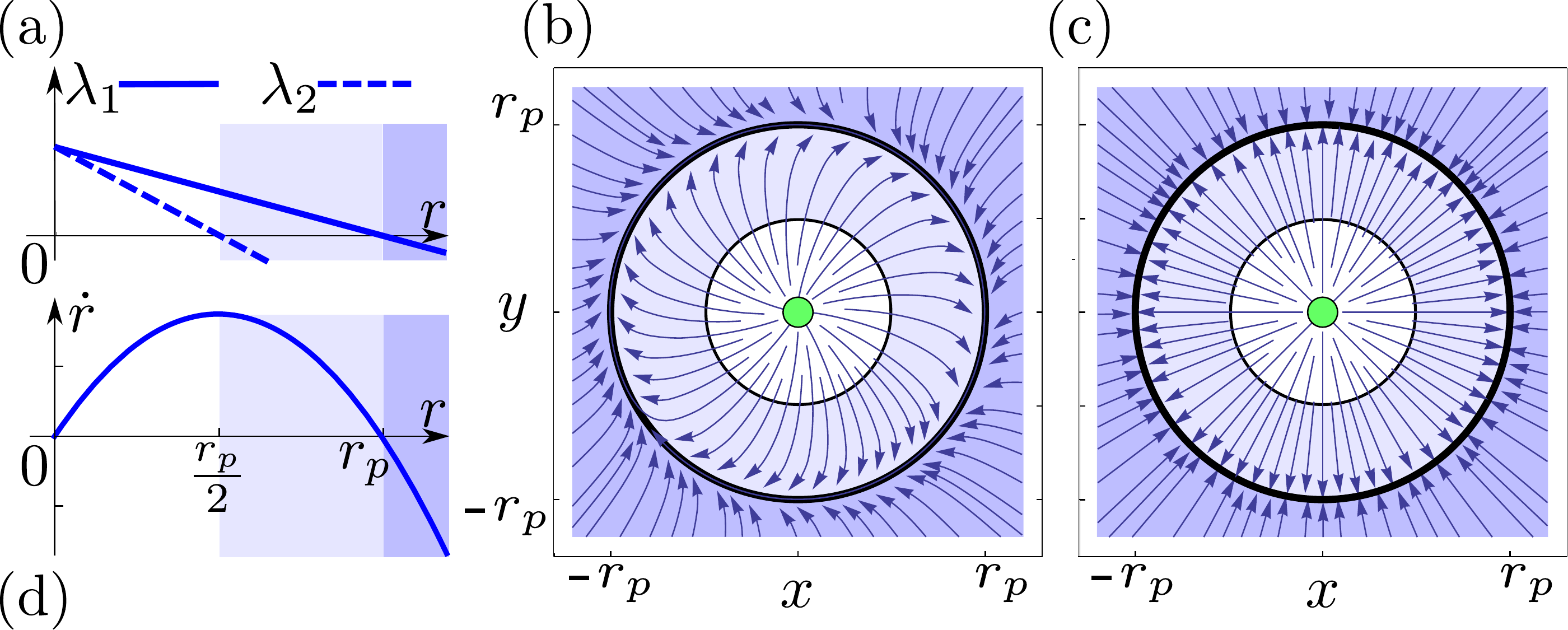}
\caption{(color online) Contraction analysis of an autonomous limit cycle (\ref{auto Poincare}).
(a) Eigenvalues $\lambda_1$ and $\lambda_2$ of the Jacobian $J$ (\ref{Jacobian}) which characterize a property of contraction. Blue (gray) area denotes a contraction region, where all trajectories converge. In the light blue (light gray) area only one of the eigenvalues is negative, and some trajectories diverge. In the white area both $\lambda_1$ and $\lambda_2$ are non-negative and all trajectories diverge.
(b) Phase portrait. Green (gray) point in the middle is an unstable node.
(c) Phase portrait in rotating reference frame.
(d) Radial velocity in (\ref{auto Poincare}) as a function of radius.
}
\label{fig:types2}
\end{figure}

When the coupling is absent $(\varepsilon_{_{A}}=0)$ the system (\ref{general example}) becomes an autonomous isochronous limit-cycle oscillator, described in polar coordinates $(r,\psi)$ as
\begin{equation}
\begin{array}{l}
 \dot{r}=\varepsilon_{_{\Gamma}}(r_{p}-r)r,\\
 \dot{\psi}=\omega_0.\\
 \end{array}
 \label{auto Poincare}
\end{equation}
Using (\ref{general example}) with $\varepsilon_{_{A}}=0$ one can find eigenvalues $\lambda_1$ and $\lambda_2$ of the Jacobian $J,$ (\ref{Jacobian}), which characterizes contraction. These eigenvalues depend only on radius $r=\sqrt{x^2+y^2},$ and are shown in Fig.~\ref{fig:types2}(a). The phase portrait of such a system is shown in Fig.~\ref{fig:types2}(b) in the laboratory coordinate system, and in Fig.~\ref{fig:types2}(c) in a reference frame rotating with an angular velocity $\omega_0.$
In the white area in Fig.~\ref{fig:types2}(a)-(d) $(\lambda_1>0,~\lambda_2>0)$ the distance between any two points can only increase in time. In the light blue (light gray) area a radial component between any 2 points will decrease in time due to the negative slope of $\dot{r}$ shown in Fig.~\ref{fig:types2}(d).
Despite the existence of a contraction region $C$ (shown as a blue (gray) area in Fig.~(\ref{fig:types2})), this system does not have a point attractor. This is because there does not exist an area $A^{\prime}\subset C$ such that states of the system do not leave it. Instead, points will eventually approach a limit cycle, which does not belong to a contraction region due to the zero eigenvalue $\lambda_1=0$ on a limit cycle.

\subsection{The case with coupling. Existence and types of chronotaxicity}
\begin{figure}[]
 \includegraphics[width=\columnwidth]{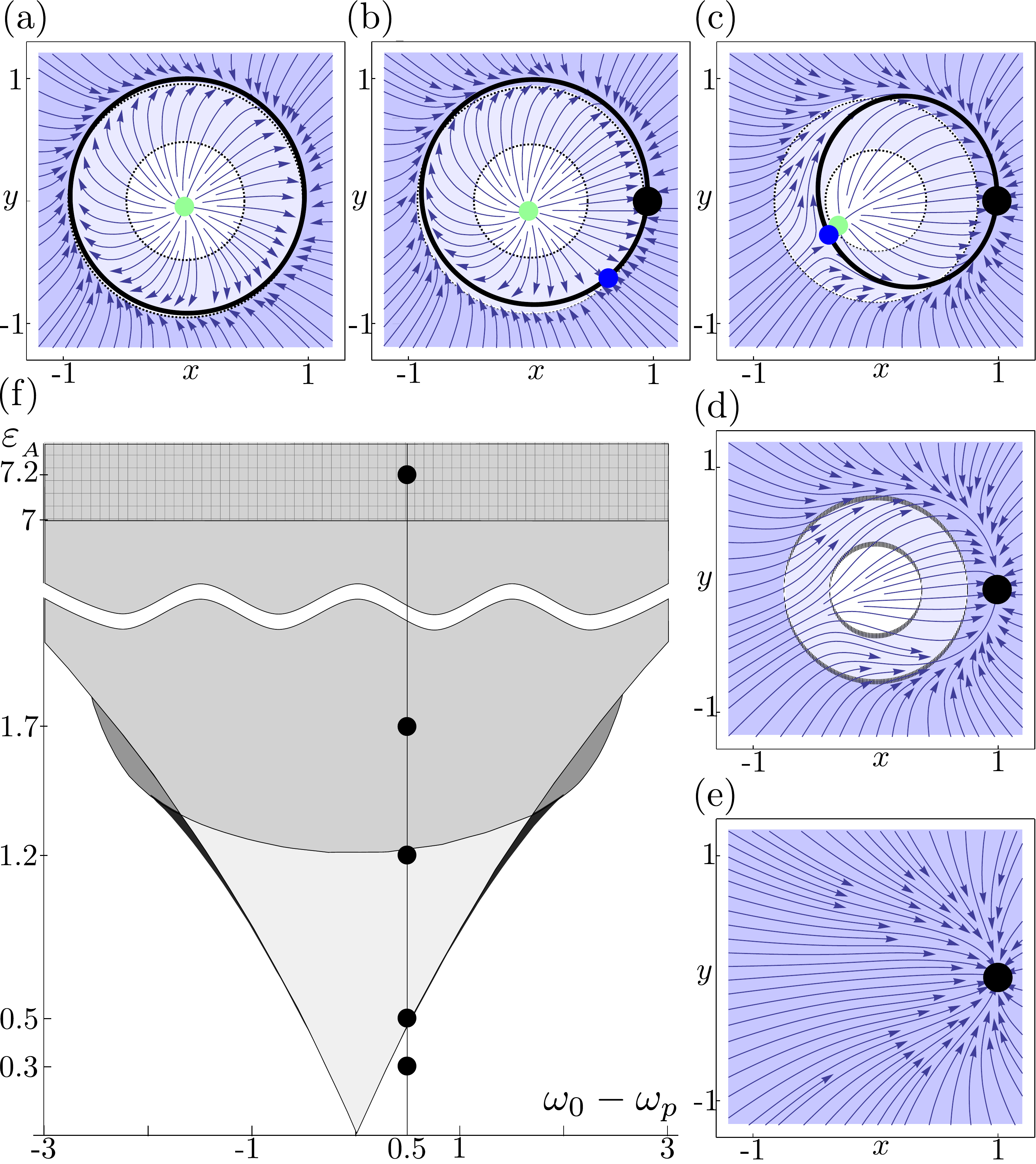}
\caption{(color online) (a)-(e) Phase portraits of the system (\ref{system polar}) with different types of chronotaxic dynamics in the rotating reference frame. The color scheme is the same as in Fig.~\ref{fig:types2}. The green (gray) point is an unstable node, the blue (black) point is a saddle point, the big black point is a point attractor. Closed black line in (a)-(c) is $\Gamma$ which is the nonautonomous analogue of a limit cycle. (f) Regions of different chronotaxic dynamics: light gray area with dots at $\varepsilon_{_{A}}=0.5$ and $1.2$ is where a point attractor and $\Gamma$ exist;  gray area with the dot at $\varepsilon_{_{A}}=1.7$ is where only a point attractor and a region of divergency exist; gray area with a small grid on top with the dot at $\varepsilon_{_{A}}=7.2$ is where only a point attractor exists. The system is not chronotaxic in the white area containing the dot at $\varepsilon_{_{A}}=0.3$. The condition of contraction is not fulfilled in the narrow black regions attached to the light gray area with dots at $\varepsilon_{_{A}}=0.5$ and $1.2$; however the system still has a point attractor and $\Gamma$ there. The condition of contraction is also unfulfilled in the dark gray regions attached to the gray area with the dot at $\varepsilon_{_{A}}=1.7$, although the system has a point attractor in those regions.   }
\label{fig: chr region}
\end{figure}

For the case $\varepsilon_{_{A}}\neq0$ in Eq.~(\ref{general example}), we first find a chronotaxic region in parameter space assuming positive and constant parameters $r_{p},$ $\varepsilon_{_{A}}$ and $\omega_{p}.$ It is convenient to rewrite Eq.~(\ref{general example}) in polar coordinates $(r,\psi)$ in a reference frame that rotates with the angular velocity $\omega_p:$
\begin{equation}
\left\{\begin{array}{l}
 \dot{r}=-\varepsilon_{_{\Gamma}}(r-r_{p})r-\varepsilon_{_{A}}(t)
 \left(r-r_{p}(t)\cos\psi\right),\\
 \dot{\psi}=\omega_0-\omega_p(t)-\varepsilon_{_{A}}(t){r_{p}(t)\over r}\sin\psi.\\
 \end{array}\right.
 \label{system polar}
\end{equation}
A point attractor corresponds to a stationary stable node in such a reference frame, i.e. it satisfies the condition of stationarity: $\dot{r}=0,$ $\dot{\psi}=0.$ The fulfilment of this requirement and the contraction property can be checked numerically. As a demonstration we choose $\omega_0-\omega_p=0.5,$ $\varepsilon_{_{\Gamma}}=7,$ $r_p=1,$ and $\varepsilon_{_{A}}$ runs from zero to higher positive values.
At small values of $\varepsilon_{_{A}}$ a limit cycle, usually denoted as $\Gamma_0,$ transforms into a closed line $\Gamma$ in the rotating reference frame shown as the black closed line in Fig.~\ref{fig: chr region}(a)-(c). This line $\Gamma$ attracts all neighboring trajectories but, in contrast to an autonomous limit cycle, this line is not a closed isolated trajectory because it moves in the laboratory coordinate system. Some part of $\Gamma$ lies inside a contraction region, and another part lies outside the contraction region.
At the stronger coupling strength $\varepsilon_{_{A}}>\varepsilon_{c1}\approx 0.467$ a saddle-node bifurcation occurs and two fixed points (black and big black dots in Fig.~\ref{fig: chr region}(b)) appear on the ``limit cycle''-like object $\Gamma$. When the stable node from this pair lies within the contraction region -- the system becomes chronotaxic; such a case is shown in Fig.~\ref{fig: chr region}(b).

At the stronger coupling strength $\varepsilon_{_{A}}$ a saddle point approaches an unstable node, Fig.~\ref{fig: chr region}(c), until the other saddle-node bifurcation occurs at $\varepsilon_{_{A}}=\varepsilon_{c2}\approx1.214$ and these two points disappear together with $\Gamma$, Fig.~\ref{fig: chr region}(d). At the same time the non-contraction region, i.e. an area where trajectories diverge from each other, becomes smaller. Only at a sufficiently strong coupling strength $\varepsilon_{_{A}}>\varepsilon_{c3}=7$ does the whole phase space become a contraction region as shown in Fig.~\ref{fig:types2}(e).
\begin{figure*}[]
 \includegraphics[width=\textwidth]{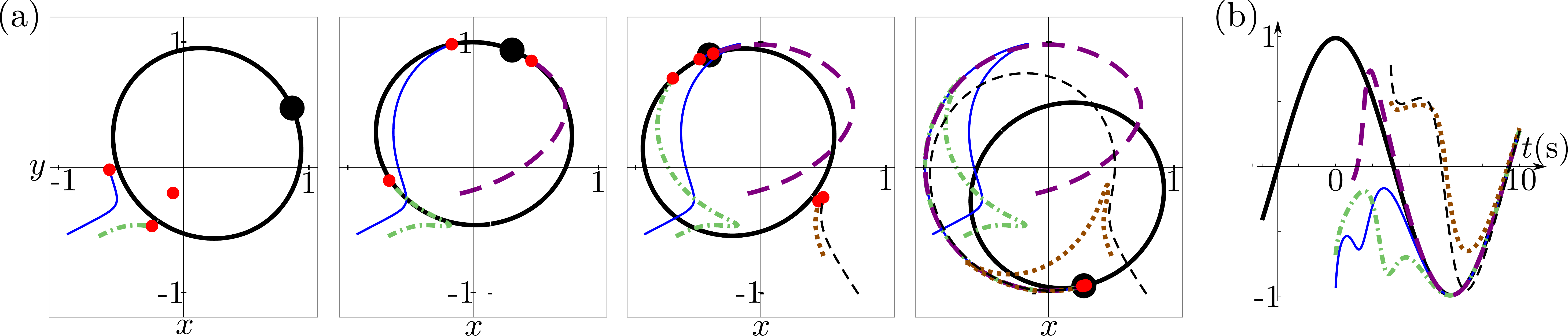}
\caption{(color online) (a) Trajectories of states indicated as red dots (dark gray) of the system (\ref{general example}) with $\varepsilon_{_{A}}=1.2,$ $\omega_0-\omega_p=0.5,$ whose phase portrait is shown in Fig.~\ref{fig: chr region}(c). Different trajectories are indicated by thin-solid line (blue), dash-dotted line (green), thick-dashed line (purple), thin-dashed line (black), and dotted (brown) line. The point attractor is shown as a big black dot and $\Gamma(t)$ -- as a thick black closed line. (b) The dynamics of the $\mathbf{x}$-component of the same system.}
\label{fig: converg of traj}
\end{figure*}

Analysis of the system at other values of $\omega_0-\omega_p$ is illustrated in Fig.~\ref{fig: chr region}(f). The colored regions (apart from narrow black and narrow dark gray areas) correspond to chronotaxic systems. The presence of $\Gamma$ and the presence of a non-contraction region distinguish 3 different types of chronotaxic dynamics: type I -- chronotaxic dynamics with a ``limit cycle''-like structure $\Gamma$ and a non-contraction region, such a dynamics being illustrated by Figs.~\ref{fig: chr region}(b)-(c); type II --  chronotaxic dynamics when a non-contraction region is present but $\Gamma$ does not exist, as shown in Fig.~\ref{fig: chr region}(d); type III -- chronotaxic dynamics where the contraction region occupies the whole space, such a dynamics being shown in Fig.~\ref{fig: chr region}(e).
There are also two regions in Fig.~\ref{fig: chr region}(f) (two narrow black and two small dark gray regions), where the condition of contraction is not fulfilled, but the point attractor exists. If transients to a point attractor are fast and negligible during observations, then the dynamics of a system in such areas can be considered as approximately chronotaxic.

Now we show how the definition of chronotaxic systems presented in this paper allows us to determine whether a system is chronotaxic when the parameters $\varepsilon_{_{A}},$ $\omega_0$ and $\omega_p$ are time-dependent functions $\varepsilon_{_{A}}(t),$ $\omega_0(t)$ and $\omega_p(t)$. First, consider the case when the value of $\omega_0-\omega_p$ is fixed and only $\varepsilon_{_{A}}$ changes. Knowing the dynamics of the contraction region demonstrated above, one then concludes that a point attractor is always inside a contraction region for any arbitrary changes in $\varepsilon_{_{A}}(t)$ inside a region of chronotaxicity.
Moreover, all transients of a system between positions of a point attractor also lie inside a contraction region. Therefore, any small vicinity of a point attractor such that it lies within a contraction region at all times can be considered as an area $A^{\prime}\subset C.$ The existence of such an area proves that the system is chronotaxic when $\varepsilon_{_{A}}(t)$ arbitrarily changes inside a region of chronotaxicity.

If only the frequency mismatch $\omega_0(t)-\omega_p(t)$ changes and $\varepsilon_{_{A}}$ is constant, then the contraction region does not change because the eigenvalues of the Jacobian (\ref{Jacobian}) do not depend on $\omega_0(t)-\omega_p(t).$ During any changes of $\omega_0(t)-\omega_p(t)$ within the chronotaxic region, a point attractor cannot leave the contraction region. One can see this by considering the following: when the value of frequency mismatch enters the region of chronotaxicity, then a point attractor appears inside the contraction region, and when the value of frequency mismatch decreases, a point attractor goes only deeper inside the contraction region. Therefore an area $A^{\prime}\subset C$ exists at all times, and the system remains chronotaxic when the value of $\omega_0(t)-\omega_p(t)$ changes arbitrarily within the chronotaxic region.

Consequently, arbitrary changes of $\omega_0(t)-\omega_p(t)$ and $\varepsilon_{_{A}}(t)$ within the chronotaxic region shown in Fig.~\ref{fig: chr region}(f) do not destroy the chronotaxicity of the system (\ref{general example}) and (\ref{system polar}). In general, one can identify a chronotaxic system from its dynamical equations. For this one should construct such an area $A^{\prime}$ which lies all the time completely within a contraction region, and the dynamical states of the system in phase space should be able only to enter this area and not to leave it.

\subsection{Response to perturbations}
The presence of an attracting line $\Gamma$ and regions where trajectories diverge from each other cause additional complexity in the dynamics, especially when external perturbations push the state of a system away from the $A^{\prime}$ region. Thus, a small instantaneous deviation may temporarily cause further deviation from a point attractor. Therefore different types of chronotaxic system will have different responses to external perturbations.

To illustrate this, consider a phase portrait for the system (\ref{general example}) with $\varepsilon_{_{A}}=1.2,$ shown in Fig.~\ref{fig: converg of traj}(a). Such a system is also shown in Fig.~\ref{fig: chr region}(c). Different initial coordinates of trajectories may be considered as resulting from instantaneous perturbations applied to a system at a point attractor. After being perturbed, all trajectories are attracted to the point attractor (big black dot), and to $\Gamma(t).$ One can see in Fig.~\ref{fig: converg of traj}(a) that the perturbed states of the system return to the point attractor in rather a complex way. When a single time-series of only the $\mathbf{x}$-component is viewed, as shown in Fig.~\ref{fig: converg of traj}(b), one can see that small deviations in the position of the system cause different dynamics and create very different trajectories, e.g. thin solid line and dash-dotted line in Fig.~\ref{fig: converg of traj}(a)-(b). At the same time, the dotted and thin dashed lines in Fig.~\ref{fig: converg of traj}(b) appear to be very similar, but in phase space they move in opposite directions along $\Gamma.$
\begin{figure}[]
 \includegraphics[width=\columnwidth]{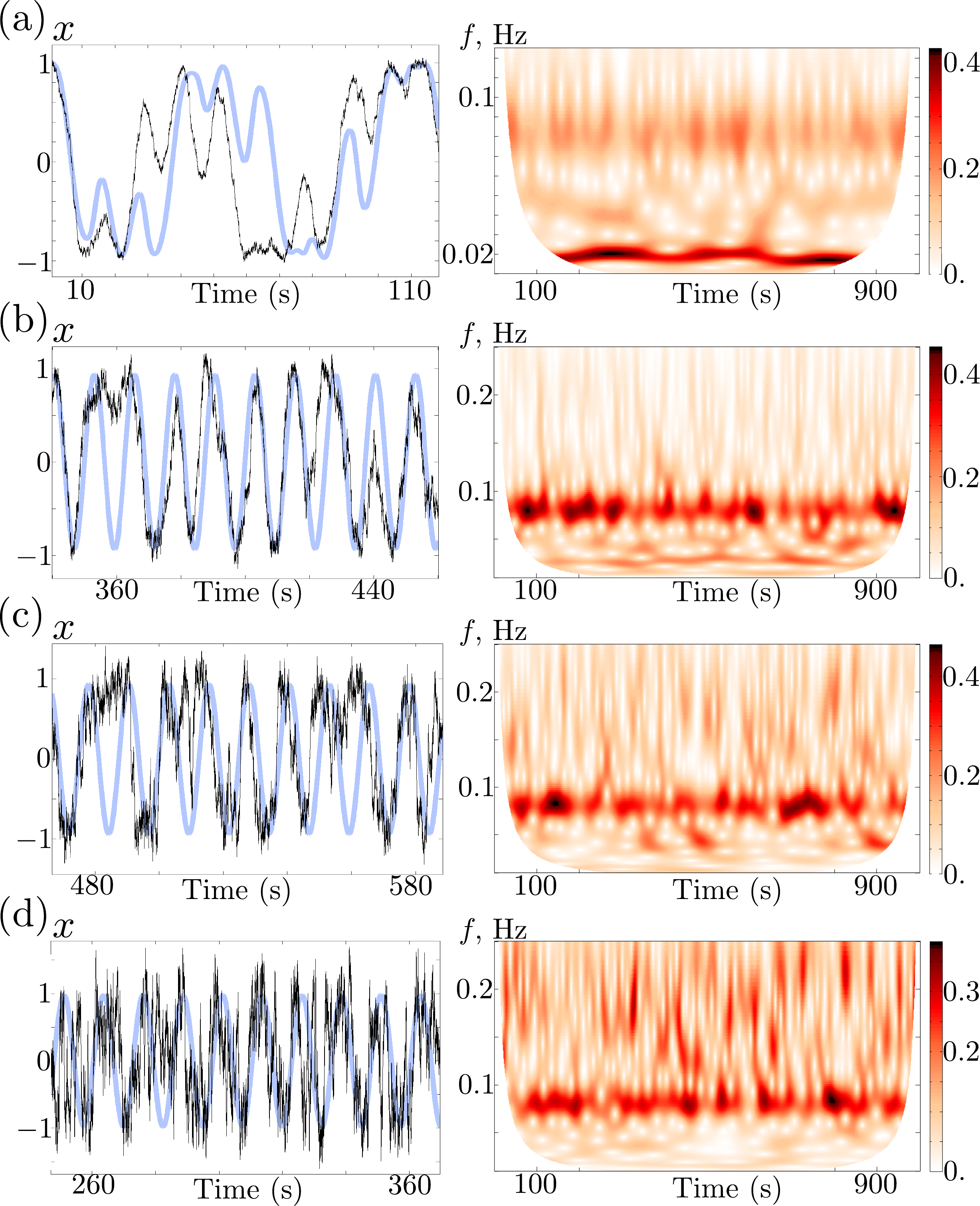}
\caption{(color online)
Response of the dynamics of (\ref{general example}) to additive white Gaussian noise (\ref{noise}). In the left-hand column: a perturbed dynamics (thin black line) and the unperturbed dynamics (thick gray line) are compared. Wavelet transforms of the perturbed dynamics are shown in the right-hand column. In all figures $f=w_p/(2\pi)\approx0.08,$  $\varepsilon_{_{\Gamma}}=7,$ $r_p=1$ and: (a) $\varepsilon_{_{A}}=0.3$, $\sigma=0.1$; (b) $\varepsilon_{_{A}}=0.47$, $\sigma=0.3$ (c) $\varepsilon_{_{A}}=0.47$, $\sigma=0.6$ (d) $\varepsilon_{_{A}}=0.9$, $\sigma=1.2.$
}
\label{fig: converg of traj2}
\end{figure}

In order to study the response to continuous perturbation, we consider the system (\ref{general example}) with additive white Gaussian noises $\xi_x(t)$ and $\xi_y(t),$ which appear on the r.h.s. of corresponding equations (\ref{general example}) for $\dot{x}$ and $\dot{y}$ respectively,
\begin{equation}
 \langle\xi_{j}(t)\rangle=0,~~\langle\xi_{i}(t)\xi_{j}(t^{\prime})\rangle=\sigma^2\delta_{ij}\delta(t-t^{\prime}),
 \label{noise}
\end{equation}
where $i,j =x,y,$ and $\delta_{ij}$ is a Kronecker delta. The dynamics perturbed by noise is shown in the left column of Fig.~\ref{fig: converg of traj2} as a black thin line, and the unperturbed dynamics (when noise is absent) is shown by a light blue (light gray) thick line. The wavelet transforms of the perturbed dynamics are shown in the right-hand column of Fig.~\ref{fig: converg of traj2}.

Figure~\ref{fig: converg of traj2}(a) shows the dynamics of a system (\ref{general example}), which does not have a point attractor (due to a weak coupling $\varepsilon_{_{A}}=0.3$). One can see that the perturbed dynamics shifts and strongly deviates from the unperturbed dynamics. The strong dominant low frequency component in the wavelet transform with $f\approx0.02~{\rm Hz}$ does not correspond to any underlying stable deterministic dynamics.

The response to noise of a chronotaxic system (\ref{general example}) with a point attractor and $\Gamma(t)$ is shown in Fig.~\ref{fig: converg of traj2}(b). In the left-hand picture around $t=360$~s one can identify the analogue of a phase slip -- the motion of a system along $\Gamma(t)$ such that a system returns to a point attractor after making a full loop along $\Gamma(t).$ It happens due to perturbations $(\sigma=0.3)$ which are weak comparatively to the attraction towards $\Gamma.$ The wavelet transform contains a strong dominant line which corresponds to the unperturbed dynamics $(f\approx 0.08~{\rm Hz})$. The dynamics of the chronotaxic system with the same coupling $\varepsilon_{_{A}}$ but in the case of a stronger perturbation $(\sigma=0.6)$ is shown in Fig.~\ref{fig: converg of traj2}(c). In such a case both the amplitude and phase are strongly perturbed. The dominant line in the wavelet transform is more hidden by the noise, but is still clearly visible.

The dynamics of a system with a stronger coupling $(\varepsilon_{_{A}}=0.9)$ and stronger perturbation $(\sigma=1.2)$ is shown in Fig.~\ref{fig: converg of traj2}(d). The time series is complex and noisy. Nevertheless, even under such perturbation it is possible to see a dominant frequency line in the wavelet transform, despite is being strongly masked by noise. The identification of the system as chronotaxic allows one to consider the dominant line in the wavelet transform as a signature of a deterministic and stable dynamics \cite{Clemson:14b}.

Thus, the dynamics of chronotaxic systems, even when it looks complex and stochastic, in principle can still allow the extraction of the deterministic component. However, if perturbations are too strong, the unperturbed dynamics becomes undetectable.

\section{Realization of chronotaxic systems}
The definition presented in this paper allows chronotaxic systems to be realized in various ways, making their practical application more probable. The drive and response systems in (\ref{1}) can be totally different dynamical systems with different dimensions.
If the drive system $\mathbf{p}$ in Eq.~(\ref{1}) is known, one can say that the $\mathbf{p}$ and $\mathbf{x}$ systems are generally synchronized. It follows from the fact that the configuration (\ref{1}) with the asymptotically stable driven system was shown to lead to general synchronization \cite{Kocarev:96}. It means that, as time goes to infinity, the states of the drive system $\mathbf{p}$ and response system $\mathbf{x},$ $(\mathbf{p}\in R^{n},~\mathbf{x}\in R^{m}),$ become connected by a static functional relationship $\mathbf{H}:R^{n}\rightarrow R^{m},$ i.e. $\mathbf{x}=\mathbf{H}(\mathbf{p}).$

However one of the simplest realizations of chronotaxic systems can be obtained when the drive system $\mathbf{p}$ is the replica of the response system $\mathbf{x}$. The corresponding example of a chronotaxic system is presented by a replica synchronization where the stability in the dynamics of a response system was also achieved, e.g. in chaotic unidirectionally coupled systems \cite{Pecora:90,Pecora:91,He:92,Pikovsky:01,Stankovski:14}. It is also important to stress that the chronotaxic system $\mathbf{x}$ is not chaotic. It is provided by the presence of a time-dependent point attractor $\mathbf{x}^{A}(t)$ and consequent insensitivity to initial conditions. Nevertheless the dynamics of the point attractor $\mathbf{x}^{A}(t)$ may also look chaotic if the drive system $\mathbf{p}$ has chaotic dynamics, for example see Ref.~\cite{Pecora:90}.

\begin{figure}[]
 \includegraphics[width=\columnwidth]{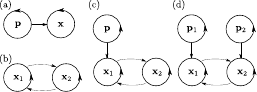}
\caption{Different configurations of interacting oscillatory systems which can lead to chronotaxic dynamics. (a) Unidirectional coupling as in Eq.~(\ref{1}) is a typical example of a chronotaxic system as discussed in Section III. (b) Bidirectional coupling. In certain cases one of the subsystems may appear chronotaxic temporarily. (c) A case where either $\mathbf{x}_1$ or both systems, $\mathbf{x}_1$ and $\mathbf{x}_2,$ can be chronotaxic and perturbed by each other. (d) Interacting chronotaxic systems. More complex yet chronotaxic dynamics is also possible.}
\label{fig: coupling}
\end{figure}

Following Eq.~(\ref{1}), a realization of a chronotaxic system $\mathbf{x}$ is shown schematically in Fig.~\ref{fig: coupling}(a). One can consider the stability of such a structure with respect to the presence of an influence of system $\mathbf{x}$ on system $\mathbf{p}.$ In that case, system $\mathbf{p}$ enters the picture equally to system $\mathbf{x},$ see Fig.~\ref{fig: coupling}(b), where it is more reasonable to use $\mathbf{x}_1$ and $\mathbf{x}_2$ instead of $\mathbf{p}$ and $\mathbf{x}.$ Thus, the dynamics will not resist external perturbations in the sense that, after the perturbations are switched off, the dynamics will not return to the unperturbed dynamics.

However, in the configuration shown in Fig.~\ref{fig: coupling}(b), the dynamics of one of the subsystems, e.g. $\mathbf{x}_2,$ may look chronotaxic temporarily; however it will not be chronotaxic rigorously. This case occurs if the influence from $\mathbf{x}_2$ to $\mathbf{x}_1$ creates only weak perturbations to the dynamics of $\mathbf{x}_1,$ which in turn creates only weak perturbations to the dynamics of the $A^{\prime}_2$ region in the phase space of $\mathbf{x}_2.$ In such a case the time interval when $\mathbf{x}_2$ looks chronotaxic should be smaller than the characteristic time of changes in $A^{\prime}_2$ due to the influence of $\mathbf{x}_2$ on $\mathbf{x}_1.$ Such a realization is expected to be especially common in real systems, where the unidirectional coupling cannot be realized. In such a situation the system may be approximated as chronotaxic during short time intervals with its structure shown in Fig.~\ref{fig: coupling}(a). This allows a time-scale separation in complex systems, and the identification of drive and response systems at certain time scales.

The chronotaxic system can be in contact with other systems. The direct influence on the drive system is not considered here, because the resulting dynamics can then be considered as a new drive system. In the configuration shown in Fig.~\ref{fig: coupling}(c) various options are possible. If the influence from $\mathbf{x}_2$ is weaker than the influence from $\mathbf{p},$ then the system $\mathbf{x}_1$ remains chronotaxic under external perturbation from $\mathbf{x}_2$. The system $\mathbf{x}_2$ can be chronotaxic too. It is possible if the influence of $\mathbf{x}_2$ on $\mathbf{x}_1$ can be suppressed by a stronger drive system $\mathbf{p}$.

Fig.~\ref{fig: coupling}(d) shows a sketch of the two interacting chronotaxic systems. Here different couplings compete with each other, so various combinations are possible. In particular, both the $\mathbf{x}_1$ and $\mathbf{x}_2$ systems can remain chronotaxic; however their interaction can make the whole dynamics look very complex. If the systems sometimes push each other out from their $A^{\prime}$ regions they still may be considered as chronotaxic provided that they spend most of their time in the $A^{\prime}$ region. In such a case when only one of these chronotaxic systems is observed, it is possible to extract its underlying deterministic dynamics which could approximate an unperturbed dynamics. It is therefore also possible to extract information about the influence, i.e. the dynamics of another chronotaxic system in this configuration.

\section{Summary}

In this work we have developed a generalized model of chronotaxic systems which are high dimensional and cannot be split into the decoupled one-dimensional dynamical systems, in contrast to previous works \cite{Suprunenko:13,Suprunenko:14}. Such systems are characterized by a time-dependent point attractor which exists in the time-dependent contraction region. Chronotaxic systems are capable of resisting to continuous external perturbations while having a complex time-dependent dynamics. Previously such systems with complex time-dependent dynamics were often treated as stochastic and the deterministic component was ignored or misidentified. The theory of chronotaxic systems, presented in this paper and in Refs.~\cite{Suprunenko:13,Suprunenko:14}, together with corresponding inverse approach methods \cite{Clemson:14b} developed to tackle chronotaxic systems, makes it possible to identify the underlying deterministic dynamics and to extract it. The resultant reduction of complexity may be useful in various applications.
The stability properties of chronotaxic systems with complex stochastic-like dynamics can help us to understand the structure and function of such systems and their interactions with the external environment. Such applications may be useful especially in living systems.
For example, features of chronotaxicity were found recently in heart rate variability in the human cardiorespiratory system with paced respiration \cite{Suprunenko:13,Clemson:14b}. As a result, the complex chronotaxic dynamics of heart rate variability was reconstructed as deterministic and an underlying deterministic dynamics was extracted. Thus, we anticipate that chronotaxic systems will find many useful applications in various research fields.

\begin{acknowledgments}
Our grateful thanks are due to P. T. Clemson, J. Gratus and P. V. E. McClintock for valuable discussions and for D. Iatsenko and T. Stankovski for valuable comments on the manuscript. This work was supported by the Engineering and Physical Sciences Research Council (UK) [Grant No. EP/100999X1].
\end{acknowledgments}

\end{document}